\documentclass[journal,onecolumn,draftclsnofoot]{IEEEtran}
\usepackage{amssymb, amsmath, bm, graphicx, color, colortbl, xfrac, tabularx}
\DeclareGraphicsExtensions{.pdf,.png,.jpg, .eps}

\begin{document}

\title{Extending the known families of scalable Huffman sequences}

\author{Timothy~C.~Petersen*,
        David~M.~Paganin,
        Imants~D.~Svalbe
  \thanks{All authors are with the School of Physics and Astronomy, Monash University, Victoria 3800, Australia.}
  \thanks{Timothy~C.~Petersen is primarily with the Monash Centre for Electron Microscopy, Monash University, Victoria 3800, Australia.}
}


\maketitle

\begin{abstract}

A canonical Huffman sequence is characterized by a zero inner-product between itself and each of its shifted copies, except at their largest relative shifts: their aperiodic auto-correlation then becomes delta-like, a single central peak surrounded by zeros, with one non-zero entry at each end. Prior work showed that the few known families of Huffman sequences (of length $N = 4n-1$, for integers $n > 1$, with continuously scalable elements) are based upon Fibonacci polynomials. Related multi-dimensional ($nD$) Huffman arrays were designed, as well as non-canonical quasi-Huffman arrays that also possess delta-like auto-correlations. We examined links between these discrete sequences and delta-correlated functions defined on the continuum, and provided simple non-iterative approaches to successfully deconvolve $nD$ data blurred by diffuse Huffman arrays.
Here we  describe new constructions for canonical Huffman sequences. Examples of length $N = 4n+1$, $N = 2n$ and families of arbitrary length are given, including scaled forms, as well as for Fibonacci-based arrays with perfect periodic auto-correlations, that are zero for all non-zero cyclic shifts. A generalization to include canonical sequences with complex scale factors invokes an equally useful dual form of delta-correlation.  We also present $1D$ arrays with a much smaller dynamic range than those where the elements are built using Fibonacci recursion. When Huffman arrays (that are comprised of inherently signed values) are employed as diffuse probe beams for image acquisition, a new two-mask de-correlating step is described here that significantly reduces the total incident radiation dose compared to a prior method that added a positive pedestal-offset.

\end{abstract}

\maketitle
\IEEEpeerreviewmaketitle

\section{Introduction}
In 1962, Huffman \cite{Huffman1962} defined what we have called the canonical form of discrete aperiodic correlation. For $1D$ finite-length Huffman sequences $H$, all off-peak auto-correlation terms are zero, except for the unavoidable single product terms at the left and right ends. For integer sequences, those end-point correlations have magnitude 1, to keep them as small as possible. In many applications, the aperiodic correlation properties of sequences are more relevant than those of the more frequently considered (and simpler) periodic case. The Huffman sequences then mirror, as closely as possible, the correlation and spectral properties of the discrete delta function: all $H$ sequences have very close to flat Fourier amplitudes. 

$H$ sequences are generated by an underlying algorithmic structure, yet their auto-correlations are more delta-like than that of finite sequences of uncorrelated white noise. Golay \cite{Golay1972,Golay1975a,Golay1975b} showed that any sequence of symmetric real-valued terms with asymmetric signs arranged in patterns like $[a,b,c,d,e,\cdots,-e,d,-c,b,-a]$ will have zeros at all even terms in its auto-correlation. For all the odd terms to also be zero requires judicious choices for the values $a,b,c,d,e,\cdots$. Hunt and Ackroyd \cite{Ackroyd1972, HuntAckroyd1980} gave several early examples of such canonical Huffman sequences, with integer values. Recent work \cite{SvalbeTCI2020} produced a large pool of scalable (with real-valued scale parameter $s$) canonical integer sequences of length $4n-1$, where the terms are constructed using the Fibonacci sequence. Extension to continuously scalable forms was made possible by generalizing the elements to Fibonacci polynomials, as implicit in the derivation by Hunt and Ackroyd \cite{HuntAckroyd1980}. An elementary proof of the auto-correlation properties for the generalized Fibonacci construction has appeared in \cite{svalbeHuffProofs}, where the Fourier transforms derived for these Huffman sequences were suitably altered to uncover several new families of odd- and even-length canonical sequences unrelated to the Fibonacci-based series. Canonical $nD$ integer arrays were also constructed in \cite{SvalbeTCI2020}, as well as what are termed `quasi-Huffmans' where their aperiodic auto-correlation closely  approximates the canonical form and retains the characteristic flat Fourier spectra.

The motivation for building quasi-Huffman sequences and arrays is to reduce the dynamic range of the $H$ sequence values. For example, the largest Fibonacci element in the sequences of Hunt and Ackroyd \cite{HuntAckroyd1980} increases significantly with the sequence length, which makes such long canonical $H$ sequences less useful. For practical applications with real signals, sequences that contain few, preferably no zeros, and that have values as constant and small as possible are strongly preferred.
The well-known binary and generalized Barker sequences \cite{BorweinJedwab2004,Jedwab} are examples of quasi-Huffman arrays. There are many related unit amplitude complex sequences such as the Chu sequences \cite{Chu1972}, often with perfect periodic auto-correlation \cite{PerfectPolyphase1992,PerfectPolyphaseGolomb1993,PerfectBlakeTirkelPolyphase}, that find practical application. However, the practical requirement for quantization \cite{ChangGolombBarkerGener1996} of the Fourier phase of each element in delta-correlated complex sequences also imposes dynamic range limitations. On the other hand, changing just the center term of the integer Barker length 13 sequence $B_{13} = [1, 1, 1, 1, 1, -1, -1, 1, 1, -1, -1, 1, -1]$ from $-1$ to $-2$ increases the auto-correlation metrics (such as the merit factor \cite{JedwabMeritFactor}) by about $25$ percent and improves the spectral flatness.  Similarly, although no $B_9$ Barker sequence exists, the sequence $[1, 1, -1, -3, -1, 1, -2, 1, -1]$ has all off peak auto-correlation entries with magnitude less than or equal to unity.

When using a $2D$ array in a real imaging system to mimic a delta-like point spread function or beam, the incident intensity of the diffuse probing beam can be minimized to reduce the impact of radiation damage and specimen distortion.  While this was the principal motivation for our Huffman diffraction physics study in Ref.~\cite{SvalbeTCI2020}, here we instead focus upon expanding the known canonical Huffman sequences, including quasi-Huffman extensions and an improved method for de-correlation.  Proofs of the delta-correlation properties for the new sequences herein are not given but various metrics of interest can be straightforwardly computed from the exact forms that are provided.

This paper is laid out as follows. Section~\ref{sec: Decorrelation} presents a method to greatly reduce the net beam intensity of $nD$ arrays by splitting the signs or Fourier phases of the incident beam. Section~\ref{sec: KroneckerHuffman} shows how the outer- or Kronecker-product of $H$ sequences, as used in \cite{SvalbeTCI2020} to produce $nD$ arrays, can also produce a large number of $1D$ quasi-Huffman sequences of close to arbitrary length. Section \ref{sec:CanonicalHuffmans} presents the algorithmic construction of several scalable $H$ sequences of length $4n+1$ and $2n$, as well as arbitrary positive-integer lengths. Perfect arrays can be constructed from Huffman sequences, as demonstrated in Sec.~\ref{sec:PeriodicHuffmans}. Section~\ref{sec:CosineHuffmans} shows that the Golay asymmetric sign pattern used to obtain sequences with all even zero term auto-correlations can be modified by adding internal patterns of reflected signs. We show here some examples of skew-symmetric $H$ sequences with long and short lengths and small amplitudes.  A dual notion of delta-correlation is described for complex scaling of Huffman sequences in Sec.~\ref{sec: ComplexHuffmanScaling}.   Section~\ref{sec: Conclusion} provides a summary and an overview of future work.

\section{De-correlation of Huffman arrays as a pair of rectified masks} \label{sec: Decorrelation}

Assume the auto-correlation of a given Huffman array $H$, or variant thereof, closely approximates a delta-function $\delta$ such that
\begin{equation}\label{CanonicalCondition}
H\star H \approx \delta,    
\end{equation}
where $\star$ represents cross-correlation. Suppose also that the blurry total measured signal $S_{T}$ comprises convolution between $H$ and a desired object signal $O$,  
\begin{equation}
S_T = O * H.    
\end{equation}
In imaging scenarios where $H$ is to be defined by an illuminated mask $M$, the diffuse probe intensity is necessarily non-negative. We previously have described an approach for de-correlating $H$ from $S_T$, involving a pair of such positive masks, by introducing a pedestal offset $\kappa$ to $H$ \cite{SvalbeTCI2020}.  Two sequence measurements using $H+ \kappa$ and $\kappa-H$ then enable de-correlation.  A possible problem with that measurement proposal can arise if $H$ has negative values with large magnitudes, as is the case for long-length Fibonacci based Huffman arrays with large scale parameters.  By offsetting the magnitudes of all element values, the required pedestal increases the implied radiation dose in the diffuse probes defined by the offset $H$ masks.

An alternative proposed here is to instead use two non-negative masks $M_1$ and $M_2$ that encode the magnitudes of the positive and negative Huffman elements, respectively, such that
\begin{equation}
H = M_1-M_2.     
\end{equation}
The total signal can then be decomposed as the difference between two successive measurements $S_1$ and $S_2$, 
\begin{equation}\label{decorrSplit}
S_T =S_1-S_2=O* M_1-O* M_2 = O* H. 
\end{equation}
The desired object signal $O$ can then be de-correlated by cross-correlating $S_T$ with the known $H$. 

For a $19 \times 19$ canonical $2D$ Huffman array of Fibonacci type \cite{SvalbeTCI2020} with unit-scaling, where the most negative value is -1764, the implied total dose of two pedestal masks is 1,273,608, as obtained by summing over the now all-positive Huffman elements.  In contrast, the split-sign masks of Eq.~\ref{decorrSplit} impart an equivalent dose of 51,076, which is about 25 times less and is equal to the sum of all elements in $M_1$ and $M_2$.  The reason for this marked improvement is that long-length canonical Huffman arrays such as $H_{19}$ contain substantial oscillations in values, which requires a large $\kappa$ to offset the negative values in the previous method.

For complex sequences or arrays, the encoded signal can be partitioned more finely, for example, into $4 \times 90$ degree phase segments. Under this view, a complex mask could be separated into four real non-negative masks, $M_{1,\textrm{Re}}, M_{2,\textrm{Re}},M_{1,\textrm{Im}}$ and $M_{2, \textrm{Im}}$, such that
\begin{equation}\label{decorrSplitComplex}
S_T =O*[M_{1,\textrm{Re}}-M_{2,\textrm{Re}}+ i(M_{1,\textrm{Im}}-M_{2,\textrm{Im}})] = O*H. 
\end{equation}
For delta-correlated complex $H$, de-correlation is achieved by cross-correlation of $H$ with $O*H$ of Eq.~\ref{decorrSplitComplex}.

\section{Varied length $1D$ delta-correlated arrays from canonical $2D$ Huffman arrays} \label{sec: KroneckerHuffman}

Multi-dimensional $nD$ canonical Huffman arrays were defined previously using outer-products of $1D$ canonical Huffman arrays \cite{SvalbeTCI2020}.  For $2D$ arrays, this procedure results in nine non-zero auto-correlation elements comprising values 1 (four times),  $-A_0$ (four times) and $A_0^2$ (once), where $A_0$ is the peak auto-correlation element of the $1D$ Huffman at zero shift \cite{SvalbeTCI2020}.  Here we report an unfolding or unwrapping of these canonical $2D$ arrays, as defined by a Kronecker product $\otimes$ instead of the outer-product $\star$.  The resulting $1D$ Huffman arrays are not technically canonical but consist of the same auto-correlation elements as the parent canonical $2D$ Huffman array and hence are highly delta-correlated. 

A large family of such Huffman arrays $H_N^K(s)$ can be generated using Kronecker products between any number of canonical Huffman arrays, where $s$ is a parameter that scales the elements of the arrays. Throughout, $s$ is taken to be real-valued, with a complex extension promoted in Sec.~\ref{sec: ComplexHuffmanScaling}.  The Kronecker product composition can be useful for creating delta-correlated arrays with only a small range of element values, such as the following example using the canonical $H_5(1)$ \cite{SvalbeTCI2020} and Fibonacci-based Huffman array $H_7(1)$:  
\begin{align}\label{HKronExample}
&H_{35}^K(1) = [1, 2, 2, -2, 1] \otimes [1, 2, 2, 0, -2, 2, -1]  \\ \nonumber
&=[1, 2, 2, 0, -2, 2, -1, 2, 4, 4, 0, -4, 4, -2, 2, 4, 4, 0, -4, 4, -2, -2, -4, -4, 0, 4, -4, 2, 1, 2, 2, 0, -2, 2, -1].
\end{align}
Using this procedure we can composite $1D$ delta-correlated arrays over a wide range of lengths through repeated Kronecker products of constituent Huffman sequences.  Combined with the $nD\rightarrow(n-1)D$ discrete-projection scheme outlined in \cite{SvalbeTCI2020} and rounding of constituent canonical Huffman array elements, one can construct long sequences of quasi-Huffman sequences with a low-dynamic range of element values, such as this $4$-bit $H_{86}$ example,
\begin{align}\label{H86}
&H_{86} = [-1,0,1,0,-1,0,2,-1,-1,-2,1,-2,1, 2,4,-2,-1,-2,-1,-5,2,4,6,-5,1,0,-3,-5,4,2,6,\\ \nonumber
&-3,-1,5,-4,1,3,-4,2, 5,-5,6,6,4,-3,0,2,-2,3,1,0,4,4,1,5,3,6,-3,-2,-3,-2,2,-6,-2,-6,\\ \nonumber
& 2,2,1,0,-4,3,1,3,0,-2,1,0,3,-1,0,-1,0,1,-1,1,-1].
\end{align}
\section{Canonical Huffman arrays of length $N=4n+1$, $N=2n$, and arbitrary lengths} \label{sec:CanonicalHuffmans}

Almost all of the canonical Huffman arrays in Ref.~\cite{SvalbeTCI2020} were of length $4n-1$, with the exception of $H_{5}(s)$. Instead of using a recursion-based approach to nullify off-peak auto-correlation, which leads to the Fibonacci sequence, here we have used computer algebra software to directly solve Diophantine equations that govern this canonical condition (akin to the direct methods of Hunt and Ackroyd \cite{HuntAckroyd1980}).  Conceding non-integer solutions and without necessarily invoking Golay skew-symmetry \cite{Golay1972}, a wide range of exact numerical and scalable forms can be determined, such as these $H_9(s)$ arrays        
\begin{align}\label{H9as}
H_{9a}(s) = &s[1/s,1,1/2(s-\sqrt{8+s^2}), -1-s/2(s-\sqrt{8+s^2}), \\ \nonumber
&-s^2/4\sqrt{8+s^2}),-1+s/2(s-\sqrt{8+s^2}),1/2(-s-\sqrt{8+s^2}),1,-1/s],
\end{align}
\begin{align}\label{H9bs}
H_{9b}(s) = &s[1/s, 1, s/2, 1/4 (4 + 2 s^2 - \sqrt{2} (4 + s^2)), 1/8 s (8 + 3 s^2 - 2 \sqrt{2} (4 + s^2)), \\ \nonumber
&1/4 (-4 - 2 s^2 + \sqrt{2} (4 + s^2)), s/2, -1, 1/s],
\end{align}
where $1/2(s-\sqrt{8+s^2})$ in Eq.~\ref{H9as} is equal to  $(s-\sqrt{8+s^2})/2$, rather than $1/(2(s-\sqrt{8+s^2}))$.  We use the former (standard) convention throughout, for all fractions. 

The following scaling relations for longer canonical Huffman arrays of length $4n+1$ were arrived at by reverse-engineering sets of exact numerical solutions, using a combination of pattern matching and judicious polynomial fitting with sparse bases,
\begin{align}\label{H13as}
&H_{13a}(s) = s[1/s,1,1/2(s+\sqrt{4+s^2}), 1/2(-4-s^2+s\sqrt{4+s^2}),-1/2(3+s^2)(s+\sqrt{4+s^2}), \\ \nonumber
&1/2(2+s^2-s\sqrt{4+s^2}(5+2s^2)),-\sqrt{4+s^2}(1+3s^2+s^4), 1/2(2+s^2+s\sqrt{4+s^2}(5+2s^2)),\\ \nonumber
&-1/2(3+s^2)(s+\sqrt{4+s^2}), 1/2(-4-s^2-s\sqrt{4+s^2}),1/2(-s+\sqrt{4+s^2}),1,-1/s]
\end{align}
\begin{align}\label{H13bs}
&H_{13b}(s) = [1, s, s^2/2, 1/16 s (8 + 3 s^2), 1/16 s^2 (8 + s^2), 1/64 s (-64 + 8 s^2 + s^4), \\ \nonumber
& 1/512 s^2 (-448 - 16 s^2 + s^4), -1/64 s (-64 + 8 s^2 + s^4),1/16 s^2 (8 + s^2), -1/16 s (8 + 3 s^2), s^2/2, -s, 1].
\end{align}
Note that all elements of $H_{13a}(s)$ asymptote toward integer values for sufficiently large real $s$, as does $H_{9a}(s)$. Similarly, both of these canonical sequences have all elements as complex Gaussian integers for imaginary $s = \pm 2i$, as does the  $H_{11}(s)$ sequence.  While these complex sequences do not adhere to Eq.~\ref{CanonicalCondition}, they exhibit an equally useful form of delta-correlation that will be clarified in Sec.~\ref{sec: ComplexHuffmanScaling}.

The following lengthy expression for the canonical $H_{17}(s)$ was arrived at in a similar manner to that of $H_{13a}(s)$ and $H_{13b}(s)$,
\begin{align}\label{H17s}
H_{17}(s) = [a(s), b(s), c(s), d(s), e(s), f(s), g(s), h(s), i(s), j(s), k(s), l(s), m(s), d(s), -c(s), b(s), -a(s)].
\end{align}
The scaling relations for the terms in Eq.~\ref{H17s} are given by
\begin{align}\label{H17sTerms}
&a(s) = 1, b(s) = s, c(s) = 1/2 s^2, d(s) = 1/4 s (4 + 2 s^2 - \sqrt{2}(4 + s^2)), \\ \nonumber
&e(s) = s^2 + 3/8 s^4 - (s^2 + 1/4 s^4) \sqrt{2}- T(s), 
f(s) = -s -1/2 s^3 + ( s +  1/4 s^3) \sqrt{2} - s T(s), \\ \nonumber
&g(s) = 1/2 s^2 - 1/2 s^2 T(s), h(s) = -s + (-s - 1/2 s^3 + (s + 1/4 s^3) \sqrt{2}) T(s), \\ \nonumber
&i(s) = (-s^2 - 3/8 s^4 + (s^2 + 1/4 s^4) \sqrt{2}) T(s), j(s) = -s + (s + 1/2 s^3 - (s + 1/4 s^3) \sqrt{2}) T(s), \\ \nonumber
&k(s) = -1/2 s^2 - 1/2 s^2 T(s), l(s) = -s - 1/2 s^3 + (s + 1/4 s^3) \sqrt{2} + s T(s), \\ \nonumber
&m(s) = -s^2 - 3/8 s^4 + (s^2 + 1/4 s^4) \sqrt{2} - T(s), \\ \nonumber
&T(s) = 1/8(512 s^2 + 480 s^4 + 160 s^6 + 17 s^8 -(256 s^2 + 320 s^4 + 112 s^6 + 12 s^8) \sqrt{2})^{1/2}.
\end{align}
Albeit more complicated, the scaling forms of the above $N=4n+1$ length arrays resemble the Binet forms of Fibonacci polynomial $N=4n-1$ canonical Huffman arrays.  Hence it may be possible that a simple recursion underlies these forms, using which an entire family of canonical arrays may be defined.  Another clue is that the twin-Huffman \cite{SvalbeTCI2020} products $H_{9a}(s)H_{9a}(-s)$ and $H_{13a}(s)H_{13a}(-s)$ each have integer values when $s$ takes on an integer value.  The $H_{9a}(s)$ and $H_{9b}(s)$ sequences differ by either starting and ending with oppositely signed-elements, or a matched-ending with +1; likewise for the length 13 sequences.  Canonical Huffman arrays of length 17 also exist for matched start and end values, for instance the $H^l_{17}$ example
\begin{align}\label{H17l}
&H^l_{17}=2[1/2, 1, 1, -1 + 2 \sqrt{2} - 2 \sqrt{2 - \sqrt{2}}, -3 + 4 \sqrt{2} -4 \sqrt{2 - \sqrt{2}}, 1 + 6 \sqrt{2} - 2 \sqrt{34 - 7 \sqrt{2}}, \\ \nonumber
&25 - 4 \sqrt{2} - 4 \sqrt{2 (10 + \sqrt{2})}, 79 + 16 \sqrt{2} - 2 \sqrt{1460 + 782 \sqrt{2}}, 145 + 48 \sqrt{2} - 8 \sqrt{394 + 223 \sqrt{2}}, -79 - 16 \sqrt{2}\\ \nonumber
&+2 \sqrt{1460 + 782 \sqrt{2}}, 25 - 4 \sqrt{2} - 4 \sqrt{2 (10 + \sqrt{2})}, -1 - 6 \sqrt{2} +2 \sqrt{34 - 7 \sqrt{2}}, -3 + 4 \sqrt{2} - 4 \sqrt{2 - \sqrt{2}}, \\ \nonumber
&1 - 2 \sqrt{2} + 2 \sqrt{2 - \sqrt{2}}, 1, -1, 1/2],
\end{align}
with a peak auto-correlation of approximately $P \approx 22.3$.

As with all Huffman arrays reported here, the dynamic range of $H^l_{17}$ can be reduced by quantizing the element values to be integers, while maintaining significant delta-like auto-correlation.  For example, rounding each element in Eq.~\ref{H17l} gives $[1, 2, 2, 1, -1, -1, 0, 1, 0, -1, 0, 1, -1, -1, 2, -2, 1]$ with $P=26$, giving just four off-peak auto-correlation elements with magnitude 2 and all others smaller.  Similarly, for scale factor $s = 3/4$ and an offset of all elements by $1/3$, the canonical Huffman $H_{17}(s)$ of Eq.~\ref{H17s} rounds to the ternary Barker sequence \cite{Jedwab,KrengelTernary} $[1, 1, 1, 0, -1, 0, 0, 0, 1, -1, 0, 1, -1, 0, 0, 1, -1]$ with merit factor $50/7$.

There are also canonical Huffman arrays of length $4n-1$ that are not defined by the generalized Fibonacci sequence yet appear somehow related.  For example, we have found a $H_{11}(s)$ sequence with the form:
\begin{align}\label{H11s}
&H_{11}(s) = s[1/s,1,1/2(s+\sqrt{5}\sqrt{4+s^2}),1/2(2+s^2+\sqrt{5} s\sqrt{4+s^2}, 7/2 s+s^3-1/2\sqrt{5}\sqrt{4+s^2},\\ \nonumber
&1+4 s^2+s^4,1/2(-7s-2 s^3-\sqrt{5} \sqrt{4+s^2}), 1/2(2+s^2-\sqrt{5} s\sqrt{4+s^2}),1/2(-s+\sqrt{5} \sqrt{4+s^2}),1,-1/s].
\end{align}

For certain scale parameters, such as $s=1$, $s=4$ and $s=11$, the canonical Huffman of Eq.~\ref{H11s} generates integer-valued elements.  For example, $H_{11}(1) = [1, 1, 3, 4, 2, 6, -7, -1, 2, 1, -1]$.  Also, the product $H_{11}(s)H_{11}(-s)$ gives integer values for all integer $s$.  The $H_{11}(s)$ sequence turns out to be a non-identical `twin' for the Fibonacci constructed sequence, $H^F_{11}(1) = [1, 2, 2, 4, 6, -1, -6, 4, -2, 2, -1]$, having exactly the same auto-correlation metrics, whilst the cross correlation between $H^F_{11}(1)$ and the new $H_{11}(1)$ is low.

Even-length quasi-Huffman arrays with a small range of integer values and flat Fourier spectra are possible, such as $[1, 2, 1, -2, 1, -1]$, $[1, 3, 4, 0, -3, 3, -2, 1]$ and $[1, -1, 0, 3, -6, 5, 5, 4]$.  Note the asymmetry in the magnitude of the end terms of the latter example.  Considering such asymmetry and relaxing element values to be non-integer permits an infinite family of scalable canonical Huffman arrays with arbitrary length to be constructed as follows,
\begin{align}\label{HArb}
&H^{\textrm{Arb}}_{N}(s) = [h^{\textrm{Arb}}_1(s),h^{\textrm{Arb}}_2(s),\cdots,h^{\textrm{Arb}}_{N-1}(s), h^{\textrm{Arb}}_{N}(s)],
\end{align}
where $h^{\textrm{Arb}}_n(s) = (-1)^{n}\sqrt{s}^{-(n-1)}$ for $1 < n < N$,  $h^{\textrm{Arb}}_1(s) = 1/(s-1)$ and $h^{\textrm{Arb}}_N(s) = (-1)^L s^{3/2-L/2}/(s-1)$. The scaling behavior of near the $s=1$ asymptote almost converts $H^{\textrm{Arb}}_{N}(s)$ into a Barker sequence of trivial structure, with all interior elements alternating between $\pm1$ and the pair of end elements diverging.  While no truly uniform Huffman sequence exists \cite{White1977}, the canonical condition is exactly maintained as $s$ tends to 1 in $H^{\textrm{Arb}}_{N}(s)$ with all but two elements approaching unit-magnitude.

Other scalable forms of finite even length can be found, such as $H^{e}_{4}(s)$
\begin{align}\label{H4Scaled}
H^{e}_{4}(s) = [1, s, 1/2 s (s + \sqrt{4 + s^2}), 1/2 (-s - \sqrt{4 + s^2})]
\end{align}
and the more complicated $H^{e}_{6}(s)$, with the form
\begin{align}\label{H6Scaled}
&H^{e}_{6}(s) =[1,s, (32 s^2 + 52 s^4 + 35 s^6 + 10 s^8 + s^{10} +W(s) (8 s + 14 s^3 + 7 s^5 + s^7)\\ \nonumber
&-\sqrt{X(s)/2} (3 s^2 (2 + s^2) + s^4 (2 + s^2) ) -(1 + s^2) (2 + s^2) \sqrt{(4 + s^2) X(s)/2}/ Z(s),\\ \nonumber
&1/4 s (-4 + Z(s)), 1/2 s (3 s + s^3 + W(s)), 1/2 (-3 s - s^3 - W(s))],
\end{align}
where
\begin{align}\label{H6Scaled}
&X(s) = (4 + s^2) (2 + s (3 + s^2) (\sqrt{4 + s^2} + s (3 + s (s + \sqrt{4 + s^2})))),  \\ \nonumber
&Y(s) = (12 s^2 + 7 s^4 + s^6 + 4 s (1 + s^2) \sqrt{4 + s^2} +s^3 (1 + s^2)\sqrt{4 + s^2}),  \\ \nonumber
&W(s) = (1 + s^2) \sqrt{4 + s^2},  Z(s) = -Y(s) + \sqrt{2} (2 + s^2) \sqrt{X(s)}.
\end{align}
Curiously,  the magnitudes (but not the signs) of all terms in $H^{\textrm{Arb}}_{N}$ match those of $H^{e}_{4}(s)$ and $H^{e}_{6}(s)$ for $s=-1/\sqrt{2}$ when $N=4$ and $N=6$ respectively.  

In prior work, upon analyzing the Fibonacci $H$ sequences of Hunt and Ackroyd \cite{HuntAckroyd1980}, we derived a new scalable family of canonical Huffman arrays of arbitrary length, which were also integer-valued for integer scalings and real-valued for real scaling. These were inferred by replacing the sine-modulation of the Fourier phase by that of a tangent function, followed by inverse Fourier transformation and truncation of elements smaller than unity \cite{svalbeHuffProofs}.  In this work we have realized that the original tangent-based Fourier spectrum can be analytically inverted to efficiently express a scalable family of real-valued canonical Huffman arrays $H^{\tan}_{N}(s)$, with arbitrary length, 
\begin{align}\label{HTan}
H^{\tan}_{N}(s) = &[s, (s^2-1)s^{n-1},\cdots,s^{-(L-3)/2}-s^{(L-3)/2},(s^2-1)s^{-m-1},\cdots],
\end{align}
where $\cdots$ signifies that the integer $n$ increases from $1$ and terminates at value $(L-3)/2$, with the opposite trend for integer $m$. Specifically, $H^{\tan}_{N}(s)$ here exactly corresponds to the inverse Fourier transform of the tangent-modulated Fourier result in \cite{svalbeHuffProofs} when the scaling parameter is transposed to $(-2 - s)/(-2 + s)$. A numerical example of length 7, $s=3$, is $H^{\tan}_{7}(3)= [3, 8, 24, -80/9, 8/27, 8/9, -1/3]$, which is also an instance of the `powers-of-three' canonical Huffman family described in \cite{svalbeHuffProofs}.   

\section{Fibonacci sequence perfect arrays from periodic canonical Huffman arrays} \label{sec:PeriodicHuffmans}

Finite sequences with zero periodic auto-correlation for all non-zero correlation shifts are called `perfect arrays' \cite{Luke1988,CekoPerfectArrays}. The canonical Huffman arrays from Ref.~\cite{SvalbeTCI2020} can be readily converted into perfect arrays by removing the $[+1,\cdots,-1]$ end elements (which necessarily define the non-zero off-peak terms in the canonical aperiodic auto-correlation) and then appending a new element with value zero.  Since cyclic permutations of elements do not change the periodic auto-correlation peak, this modified Huffman array can be brought into closer correspondence with the ordered Fibonacci sequence,  
\begin{align}\label{PeriodicHuffman}
H^{p}_{N-1}(s) =2 s[&X(s),F_{-M}(s),F_{-(M-1)}(s),F_{-(M-2)}(s) \cdots,F_{M-2}(s),F_{M-1}(s),F_{M}(s)], 
\end{align}
where $X(s)=F_{M+1}(s)/2 - F_{M}(s)/s$ and $M = (N-3)/2$ must be an integer. Were it not for the leading term $X(s)$, one could argue that Fibonacci sequences of length $N = 4n+3$ are perfect arrays. An example is $1/2 H^{p}_{10}(1) = [0,-3,2,-1,1,0,1,1,2,3]$, which derives from the (halved) canonical $H_{11}(1)$ Huffman sequence of~\cite{SvalbeTCI2020} with unit scale parameter.

It is worth noting that the same procedure for converting the Fibonacci canonical Huffman arrays into perfect arrays also works for the infinite family of arbitrary-length canonical Huffman arrays described earlier in Sec.~\ref{sec:CanonicalHuffmans}.  Specifically, the first and last terms in Eq.~\ref{HPerfectArb} were deleted to suppress the off-peak auto-correlation values and a new but similar element was appended to then set all of these residual values to zero,
\begin{align}\label{HPerfectArb}
&H^{\textrm{Arb}*}_{N}(s) = [h^{\textrm{Arb}*}_1(s),h^{\textrm{Arb}}_2(s),\cdots,h^{\textrm{Arb}}_{N-1}(s)],
\end{align}
where $h^{\textrm{Arb}*}_1(s) = (s-1)^{-1} (1+(-1)^{1+N} s^{1-L/2})$.

%

\section{A family of non-canonical delta-correlated arrays of length $N=4n+1$} \label{sec:CosineHuffmans}

By flipping half the signs of the Fibonacci based canonical Huffman arrays, an infinite family of delta-correlated arrays are generated with length $4n+1$.  Consisting solely of Fibonacci polynomials, these arrays $H_N^+(s)$ that start and end with $+1$ are given by,  
\begin{align}\label{HCosine}
H_N^+(s) = 2 s[(2s)^{-1},F_{1}(s),F_{2}(s),\cdots,  F_{M}(s),F_{M+1}(s)/2 - F_{M}(s)/s,-F_{-M}, \cdots,-F_{-2}(s),-F_{-1}(s),(2 s)^{-1}], 
\end{align}
where $M = (N-3)/2$. This sequence can also be derived by replacing the sine functions in the smooth Fourier-phase \cite{svalbeHuffProofs} with a cosine function.

These arrays are not canonical but, like the Kronecker derived arrays of Sec.~\ref{sec: KroneckerHuffman}, all but five elements in the auto-correlation are zero, which range over the set of elements ${1,-2\sqrt{P-2},P}$, where $P$ is the peak.  An example is $H_9^+(1)=[1, 2, 2, 4, -1, -4, 2, -2, 1]$, which has the aperiodic auto-correlation:
\begin{align}\label{HCosine}
H_9^+(1)\star H_9^+(1) =[1, 0, 0, 0, -14, 0, 0, 0, 51, 0, 0, 0, -14, 0, 0, 0, 1].
\end{align}
This pattern of negative elements at $1/4$ and $3/4$ of the auto-correlation length is generic to all $H_N^+(s)$ arrays.

\section{Dual of delta-correlation to harness complex-scaled canonical Huffman arrays} \label{sec: ComplexHuffmanScaling}
The de-correlation of complex masks was described in Sec.~\ref{sec: Decorrelation}.  While there may exist a broad class of general complex Huffman sequences, outside the scope of this study, a simple strategy is to try complex scale factors for the existing real-valued canonical Huffman sequences.  Unfortunately this idea spoils the canonical Huffman condition, since the general cross-correlation between complex sequences $f$ and $g$ is defined by a product of $\bar{f}$ and shifts of the sequence $g$, where the over-bar denotes complex conjugation. 
One remedy to restore the notion of delta-correlation is to define a dual form of auto-correlation $\bar{f}\star f$ and seek canonical Huffman arrays by analogy with Eq.~\ref{CanonicalCondition},
\begin{equation}\label{DualCanonicalCondition}
\bar{H}\star H \approx \delta.
\end{equation}
Since the altered $\bar{H}(s)$ term reverses the complex conjugation inherent in the conventional auto-correlation, all of the scaled canonical Huffman sequences here and in \cite{SvalbeTCI2020, svalbeHuffProofs} trivially remain delta-correlated in the sense of Eq.~\ref{DualCanonicalCondition}, with zero off-peak correlation elements for all but the very ends, when extended to complex scale factors.  The resulting dual-canonical Huffman sequences $H_d(s)$ do not have flat Fourier spectra but, rather, the cross-spectrum or product between the conjugate of $\mathcal{F}[\bar{H}_d(s)]$ and $\mathcal{F}[H_d(s)]$ is flat, for Fourier transform $\mathcal{F}$.

For intended applications of complex Huffman arrays such as in digital watermarking or the de-correlation of blurred data from diffuse beams \cite{SvalbeTCI2020}, the dual form of auto-correlation is just as useful as the standard canonical condition.  For example, Eq.~\ref{decorrSplitComplex} with $H_d(s)$ need only be cross-correlated with $\bar{H}_d(s)$, rather than $H_d(s)$, to suitably enact de-blurring.  Moreover this dual-form of delta correlation extends to $nD$, upon implementing outer-products of $H_d(s)$ in the same manner as in \cite{SvalbeTCI2020}. 

This extension to complex scalings offers interesting analogies with generalized Barker sequences \cite{Jedwab}.  For example, many of the canonical Huffman sequences here and in \cite{SvalbeTCI2020, svalbeHuffProofs} largely comprise elements of constant modulus when the scale parameter is set to $s=\sqrt{-1} = i$.  For instance, the arbitrary length sequence $H_{\textrm{int}}^N(s)$ from \cite{svalbeHuffProofs}, when halved and $s$ is set to $i$, has values $\pm 1$ or $\pm i$ for all but the first and last elements, which have magnitude $1/2$. An example is $H_{\textrm{int}}^7(i)/2= [i/2, -1, -i, 1, i, -1, -i/2]$.  In fact every element of $H_{\textrm{int}}^N(s)$ has unit-modulus when $s = e^{\pm \pi i 5/6}$, for any length $N$.  An example for $N=7$ is 
\begin{equation}\label{UnitModulusHuffman}
H_{\textrm{int}}^7(s=e^{\pi i 5/6})=1/2[i - \sqrt{3}, -1 - i \sqrt{3}, 
 i + \sqrt{3}, -2, -i + \sqrt{3}, -1 + i \sqrt{3}, -i - \sqrt{3}].
\end{equation}
Similarly, $H^{e}_{4}(s)$ satisfies Eq.~\ref{DualCanonicalCondition} and has unit modulus when $s = e^{\pi i/2}$, as does $H^{\textrm{Arb}}_{N}(s)$ for $s = e^{\pi i/3}$ for any length $N$.  On the contrary, while dual-correlated in the sense of Eq.~\ref{DualCanonicalCondition}, there are no complex scalings for which  the Fibonacci based canonical Huffman arrays can be made uni-modular.  Also,  the dual-correlated $H^{\tan}_{N}(s)$ has unit magnitude when $s=e^{\pi i/6}$ only for $N = 1+12 n$, due to modulus of the middle term which differs from unity for other lengths.  One could consider more general scalings of the real canonical Huffman arrays, like the quaternions, however issues such as non-commutation render this extension non-trivial and outside the scope of this work. Indeed there exist interesting constructions, including exhaustive searches \cite{QuaternionBlakePhD}, of uni-modular perfect sequences over the quaternions \cite{QuaternionKuznetsovPhD,QuaternionKuznetsovConf,QuaternionKuznetsovPaperI}, which similarly call for generalized notions of right- and left-auto-correlation which respect the absence of commutation \cite{QuaternionKuznetsovPaperII}.  Recently Bright et al. \cite{QuaternionPolyphase2020} derived new constructions for arbitrarily long perfect sequences over the quaternions, contradicting a prior existence conjecture concerning the allowable lengths given constrained alphabets of sequences in this context \cite{PerfectBlakeTirkelPolyphase}.

In short, if we ignore the conjugation operation in the definition of the cross-correlation, then one can devise the equivalent of generalized Barker sequences \cite{Jedwab} that are also canonical Huffman sequences in this sense, of arbitrary length.  Larger complex scale factors for our canonical Huffman sequences increase the magnitude of the $\bar{H}_d(s)\star H_d(s)$ peak, similar to the case of real $s$.

\section{Conclusion} \label{sec: Conclusion}

For applications where Huffman sequences must be defined by non-negative masks, such as scanned diffuse probes for imaging \cite{SvalbeTCI2020}, offsetting all elements to remove any negative values increases the implied irradiation dose.  Here we presented a useful alternative by proposing sequential scans with composite positive-valued masks. It was shown that de-correlation can then be achieved through the difference of the two encoded beam image intensities, or for four masks with complex Huffman encoding, and that this partitioning greatly reduces the implied dose. We next extended the known set of scalable canonical Huffman sequences to include examples for several $4n+1$ cases, some even length sequences, and two families of arbitrary length. We converted canonical sequences into quasi-Huffman arrays through variations of Hunt and Ackroyd's \cite{HuntAckroyd1980} Fibonacci based sequences, as well as by composing Kronecker products to create long delta-correlated sequences from constituent canonical forms. We have further shown that canonical aperiodic Huffman sequences can be converted to perfect periodic sequences.  Lastly, we altered the defining auto-correlation sum to preserve the notion of delta-correlation for 
complex scalings of all of our canonical Huffman sequences and also demonstrated that arbitrary length sequences with all elements of unit-modulus are possible under this scheme. Future work will present methods to construct families of Huffman arrays from higher order correlations.

\bibliographystyle{IEEEtran}

\bibliography{refs}
\vskip -2\baselineskip
\begin{IEEEbiographynophoto}{Imants D. Svalbe}
Imants Svalbe completed a PhD in experimental nuclear physics at Melbourne University in 1979. His current work applies Mojette and Finite Radon transforms to design $nD$ geometric structures of signed integers that act as zero-sum projection ghosts in discrete tomography and to build large families of $nD$ integer arrays that have optimal correlation properties.
\end{IEEEbiographynophoto}
\vskip -2\baselineskip
\begin{IEEEbiographynophoto}{David M. Paganin}
David Paganin received his PhD in optical physics from Melbourne University in 1999, and has been with Monash University since 2002. His research interests include x-ray optics, visible-light optics, electron diffraction, neutron optics and non-linear quantum fields. 
\end{IEEEbiographynophoto}
\vskip -2\baselineskip
\begin{IEEEbiographynophoto}{Timothy C. Petersen}
Timothy Petersen completed a PhD in the condensed matter physics of disordered carbon, at RMIT University in 2004.  He has performed experiments across a range of microscopy techniques to study disordered solids and develop new diffraction physics theories.
\end{IEEEbiographynophoto}

\end{document}